# Development of Multi-physics Finite Element model to investigate Electromagnetic Forming and Simultaneous Multi-point Perforations of Aluminium tube

Avinash Chetry* and Arup Nandy

*Department of Mechanical Engineering, Indian Institute of Technology, Guwahati, Assam, 781039, India*

*avinash.chetry@iitg.ac.in

Electromagnetic forming and perforations (EMFP) are complex and innovative high-strain rate processes that involve electromagnetic-mechanical interactions for simultaneous metal forming and perforations. Instead of spending costly resources on repetitive experimental work, a properly designed numerical model can be effectively used for detailed analysis and characterization of the complex process. A coupled finite element (FE) model is considered for analyzing the multi-physics of the EMFP because of its robustness and improved accuracy. In this work, a detailed understanding of the process has been achieved by numerically simulating forming and perforations of Al6061-T6 tube for 12 holes and 36 holes with two different punches, i.e., pointed and concave punches using Ls-Dyna™ software. In order to shed light on EMFP physics, a comparison between experimental data and the formulated numerical simulation has been carried out to compare the average hole diameter and the number of perforated holes, for different types of punches and a range of discharge energies. The simulated results show acceptable agreement with experimental studies, with maximum deviations being less than or equal to 6%, which clearly illustrates the efficacy and capability of the developed coupled Multi-physics FE model.

*Corresponding Author, E-mail: avinash.chetry@iitg.ac.in

# 1. Introduction

Perforated tubes are utilized in many practical applications to carry out certain specified design functions. For instance, perforated tubes are used in the lubrication system to support the filtration media, simultaneously holding it in place and letting the fluid flow. In addition, the oil and gas production industries use perforated tubes to attain some filtration requirement. In acoustic mufflers perforated tubes reduce generated noise, which finds application in vehicles, static engines and generators. For the pipes carrying hot fluids, one layer of outer perforated tube serves as a shield against direct contact with the hot pipe. In aerospace industries, perforated tubes are used for exhaust systems and in space-shuttle. Conventional forming and shearing of for perforated tubes lead to burrs and slivers in the shearing region [1]. Also, simultaneous forming and perforation lead to higher productivity. Thus, in recent times various research activities are happening in the field of Electromagnetic Forming and Perforation (EMFP).

Some early contributions on this subject were reviewed in details [2] including various experimental and numerical investigation of flat sheets and tubular components to understand the feasibility of the electromagnetic shearing process. Authors claimed that the process is advantageous with contact-free force applications, high flexibility for tool coil, and higher production rate. Later, in the early 21$^{st}$ century, the experimental and FE analysis to understand the blank separations mechanism due to the high strain rate process based on electrohydraulic trimming [3] and Punchless electromagnetic shearing [4] were carried out by researchers. The experimental feasibility and impact of multiple process variables on the manufacturing of the perforated tubes using electromagnetic forming and multi-point perforation using a pointed and concave punch, was carried out in [5]. Based on the experimental observations, the authors claimed that the disadvantages of the conventional perforation procedure, like formations of slivers and burr, were overcome using EMFP.

3D FEM analysis of multi-point EMFP of tubes has not received considerable attention from the research community. The complexity of the physical events involved makes complete knowledge of physics exceedingly difficult: like wave distributions of electromagnetic field, Lorentz force distribution, heating involved during plastic deformation, inertial effects during high strain rate forming, etc. [6]. Hence the development of such a multi-physics simulations model is a handy tool that provides precise explanations of the entire process. Developing of such numerical tool reduces the cost of repetitive experiments [7]. Also, such tool can be repeatedly used in different optimization or ANN models [8,9] for determining optimum geometrical and process parameters. It is reported from open literature [10-18], that the coupled and non-coupled approaches are used to understand the magnetic field intensity, Lorentz force distributions, plastic deformations behavior, shearing mechanism, etc., which are difficult to capture during the experiments because of the high-speed EMF process. The various electromagnetic forming (EMF) applications solved using coupled numerical model are Electromagnetic expansions of Al sheets [10], EM cladding of aluminium alloy to mild steel [11], Electromagnetic Crimping (EMC) of Al-SS tube to rod joint [12], Electromagnetic Crimping (EMC) of Cu-SS tube to tube joint [13], Magnetic Pulse Welding (MPW) of Steel-Al tubular parts [14],etc., and using non-coupled numerical model are electromagnetic punch-less shearing of Al tube [15], Electromagnetic Forming (EMF) of sheets [16], Electromagnetic forming of muffler tube [17], etc. A comparative study [18] states that the coupled approach using the Ls-Dyna™ software package is more accurate than non-coupled simulation strategies, but it is computationally expensive and the simulation interface is less user friendly.

Johnson cook (JC) material model [19] is the best choice for the constitutive model for pulsed electromagnetic forming since it illustrates the combined effect of strain, strain rate and temperature. However, another popular constitutive model for high speed forming is Cowper Symonds material model [20, 21] which does not consider the effect of temperature.

The present work deals with fulfilling two objectives. First, developing a Multiphysics-coupled numerical simulations model using a JC material model for forming and perforations of Al6061-T6 tube for 12 holes and 36 holes with two different punches, i.e., pointed and concave punches which will shed light into deformations behavior and shearing mechanism. Secondly, a comparative study is carried out between the developed 3D FEM model and experimental results of the average diameter of holes and the number of perforated holes for three different energy levels to illustrate the efficacy of the coupled Multi-physics FE model. The proposed FE model will provide a foundation to optimize the various process and geometrical parameters without spending much on costly resources through experiments.

**2. Working Principle and Process Analysis of EMFP**

The general setup of the multi-point EMFP of tubes is illustrated in Figure 1. When the switch is on, the energy is briefly held in the capacitor before being quickly discharged into the coil. The workpiece (conductive tube) encompasses a strong magnetic field created by the current carrying coil. This magnetic field produces eddy currents within the workpiece. The secondary magnetic field due to this eddy current opposes the primary field. The repulsive forces developed are known as Lorentz forces which act on the entire volume of the conductive tube resulting in plastic deformation at high strain rates. The rigid punches fixed in the cast-iron die surrounding the workpiece, perforate the expanding tube.

The usual arrangement for an electromagnetic forming and perforation circuit is an RLC circuit. A resonant circuit is used to mimic the electromagnetic forming machine setup which includes a capacitor $C$, an inductance $L_i$, and a resistor $R_i$. The resistance $R_{\text{coil}}$ and inductance $L_{\text{coil}}$ of the tool coil are linked in series to the EMF forming machine, and this circuit is designated as the primary circuit. The secondary circuit representing workpiece inductance and resistance $L_w$ and $R_w$ respectively, are linked to the primary circuit by mutual coupling effect as depicted in Figure 2(a). Thus, a simplified circuit diagram is obtained using the aggregate values of the workpiece system with coil designated as equivalent inductance $L_{Eq}$ and resistance $R_{Eq}$, depicted in Figure 2 (b).

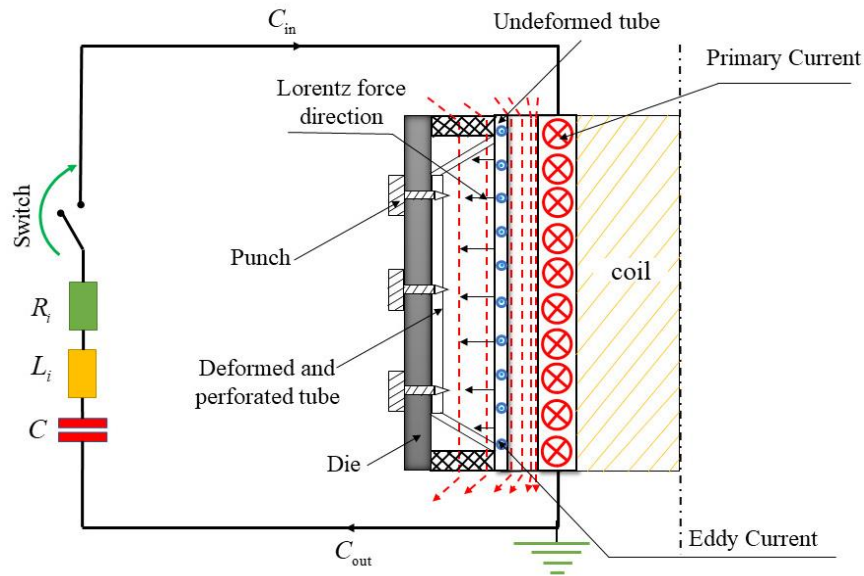

**Figure 1.** A sectional view of the setup about the axis.

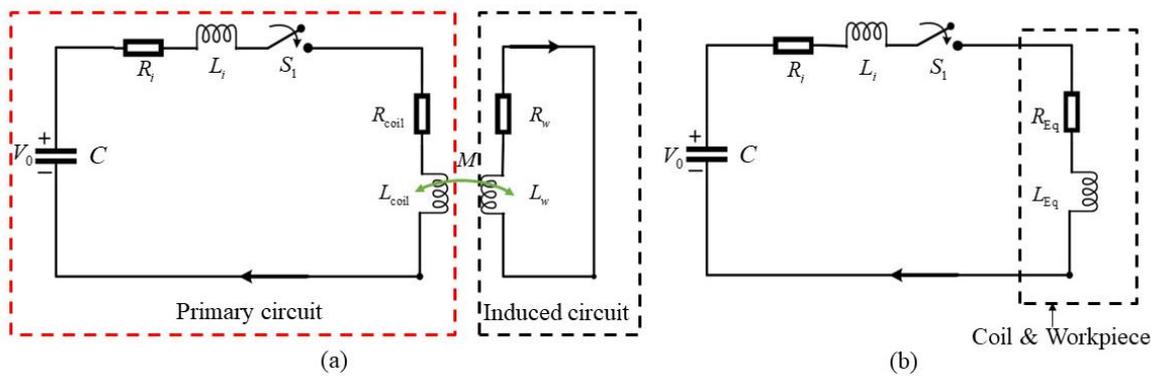

**Figure 2.** (a) Actual elementary circuit diagram of EMFP, (b) simplified elementary circuit diagram

The governing equations for discharge process in EMF while ignoring the changing of circuit parameters throughout the EMF process, given by [22],

$$L\frac{\mathrm{d}^2 I}{\mathrm{d}t^2} + R\frac{\mathrm{d}I}{\mathrm{d}t} + \frac{1}{C}I = 0, \qquad (1)$$

$$I\big|_{t=0} = 0, \qquad (2)$$

$$\left.\frac{dI}{dt}\right|_{t=0} = \frac{V_0}{L}, \tag{3}$$

$$L = L_{Eq} + L_i, \tag{4}$$

$$R = R_{Eq} + R_i, \tag{5}$$

where $I$ represent the discharge current, $V_0$ represents the voltage of the first discharge, $L$ and $R$ are the total circuit inductance and resistance, respectively. On solving the above equations, the discharge current $I$ in the simplified RLC circuit can be obtained in the form of a damped sinusoidal function given by Eq. (6) as expressed in [22],

$$I = \frac{V_0}{\sqrt{1-\zeta^2}} \sqrt{C/L} \exp(-\zeta \omega_n t) \sin(2\pi f t), \tag{6}$$

where $f$ denotes the discharge circuit's frequency as

$$f = \frac{\omega_n}{2\pi} \sqrt{1-\zeta^2}, \tag{7}$$

$\omega_n$ denotes the discharge circuit's natural frequency as

$$\omega_n = 1/\sqrt{LC}, \tag{8}$$

$\zeta$ denotes the discharge circuit's damping ratio as

$$\zeta = R\sqrt{C/L}/2. \tag{9}$$

The value of skin depth $\delta_s$ can be evaluated using electrical conductivity $\kappa$, the discharge frequency $\omega$, and permeability $\mu$ of the workpiece given as Eq. (10),

$$\delta_s = \sqrt{\frac{2}{\omega \mu \kappa}}, \tag{10}$$

which is a crucial parameter while considering the magnetic penetration assumption [2]. Current interact with magnetic field to produce the Lorentz force as

$$\vec{F} = \vec{J} \times \vec{B} = \left(\vec{\nabla} \times \frac{\vec{B}}{\mu_m}\right) \times \vec{B}. \tag{11}$$

Integrating radial Lorentz force across the tube thickness, we can have the transient magnetic pressure [2] as

$$p_{\text{m}}(r,t) = \int_{r_i}^{r_o} F(r,t)\mathrm{d}r = \frac{1}{2}\mu\left(H_{\text{gap}}^2(t) - H_{\text{pen}}^2(t)\right), \tag{12}$$

where $H_{\text{gap}}$ is the concentrated magnetic field in the radial gap among the workpiece and the tool, $H_{\text{pen}}$ is the penetrating magnetic field. If the workpiece wall thickness is greater than the skin depth, it is conceivable that the penetrating magnetic field is wholly guarded within the workpiece. Therefore, one can neglect the effect of penetrating magnetic field in Eq. (12). Then, we have

$$p_{\text{m}}(r,t) = \frac{1}{2}\mu H_{\text{gap}}^2(t). \tag{13}$$

## 3. Coupled Simulations Methodology:

The multi-physics software, Ls Dyna™ [23] is used for coupled electromagnetic and mechanical simulation. In Figure 3, a flow chart pictorially represents different steps in coupled simulation. The process is governed by two loops. In the first loop, at each time step, the EM module permits the current waveform to pass through the coil and evaluates the associated magnetic fields, magnetic forces, etc., using Maxwell equations. The calculated Lorentz force values are considered in the governing differential equations for the structural module to account for tube deformation and update the geometry, which again has an impact on the electromagnetic field. This altered EM field also has an effect on structural deformation. These iterations continue in the first loop until the convergence is reached. The second loop over time continues till the end time, updating the current time step with corresponding current value in the first loop.

The EM solver applies finite element method (FEM) to determine the EM fields inside of the conductors, whereas the boundary elements method (BEM) is utilised to handle the interaction between conductors. Depending on the size of the conductor's surface meshes, the EM solver will independently generate those boundary elements. This implies that no air mesh is required, which facilitates conductor mobility without considering air mesh deformation. One crucial aspect to note is the periodic recalculation of the FEM-BEM solver, which, by default, are only configured once during the initialization phase. The simulation uses 5000 cycles for both the FEM and BEM solvers to balance accuracy and speed. The BEM solver uses a relative tolerance of $\varepsilon_{BEM} = 1 \times 10^{-6}$, while the FEM solver uses a relative tolerance of $\varepsilon_{FEM} = 1 \times 10^{-3}$ (see Fig. 3).

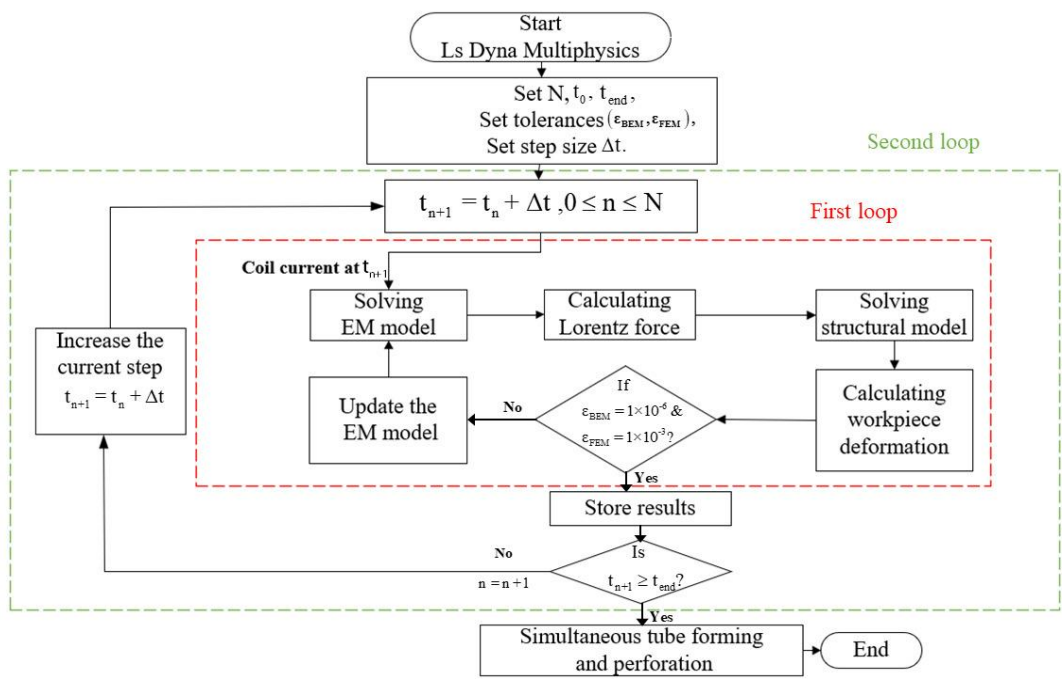

**Figure 3.** Flow chart of the coupled finite element (FE) model.

Figure 4(a), shows the 3D model with all the components, including the Al 6061-T6 tube, the Cu coil, and the SS304 punch and die, together with their respective geometrical dimensions according to the experimental studies [5]. As shown in Figure 4(b), the punches are evenly distributed in the circumferential direction throughout the tube. In the instance of a 12-punch configuration, three set of punches are placed with a mutual central distance of 20 mm. In each set, there are four punches 90° apart along the circumferential direction. In the 36-punch arrangement, circumferential spacing is changed from 90° to 30° with 12 punches in one

particular axial location. The meshing of the geometries is carried out in HyperMesh [24], and then they are imported into the Ls Dyna™ software for analysis. An eight-node linear brick element (C3D8R) fine mesh is created for aluminium tubes, and a four-node rigid element coarse mesh (R3D4) is created for coil, punch, and die. A mesh convergence study was carried out in [25], where a quarter model with single pointed and concave punch geometry was analysed to determine the optimum mesh size. Based on this, the element sizes for the analysis has been taken as 0.1 mm for the most critical components (Al6061-T6 tube) in case of 12 multi-point punch geometries and for 36 multi-point punch geometries respectively. Detailed mesh descriptions for those geometries have been presented in Table 1. It is important to carefully consider the boundary conditions and ensure that they are appropriately applied to reflect the real-world conditions of the problem being analysed. As shown in Fig. 4(c) two ends of the tube are constrained to move in x and y direction whereas axial translation (z direction) as well as all the rotations are free.

**Table 1.** Detailed of the FEM meshes.

| Arrangement | Tube element sizes (mm) | Total number of elements | Total number of nodes |
| --- | --- | --- | --- |
| 12 multi-point punch geometry | 0.1 | 43919 | 70298 |
| 36 multi-point punch geometry | 0.1 | 46295 | 93401 |

The solenoid copper coil can generate a magnetic field when an electrical current is passed through it. During simulations, the solenoid copper coil typically has two terminals or segments where the current can flow in and out, as shown in Figure 4(c). During electromagnetic forming and perforation, the tube and the punch interact with each other, and the friction coefficient is pivotal in this interaction. In this case, the static and dynamic coefficients of friction are both 0.30. In some simulations, Joule heating loss is not considered due to the high strain rate process [13, 17, 25]. Similarly, in some simulations of electromagnetic forming and perforation, isotropic conductance of materials is assumed [13, 17, 25].

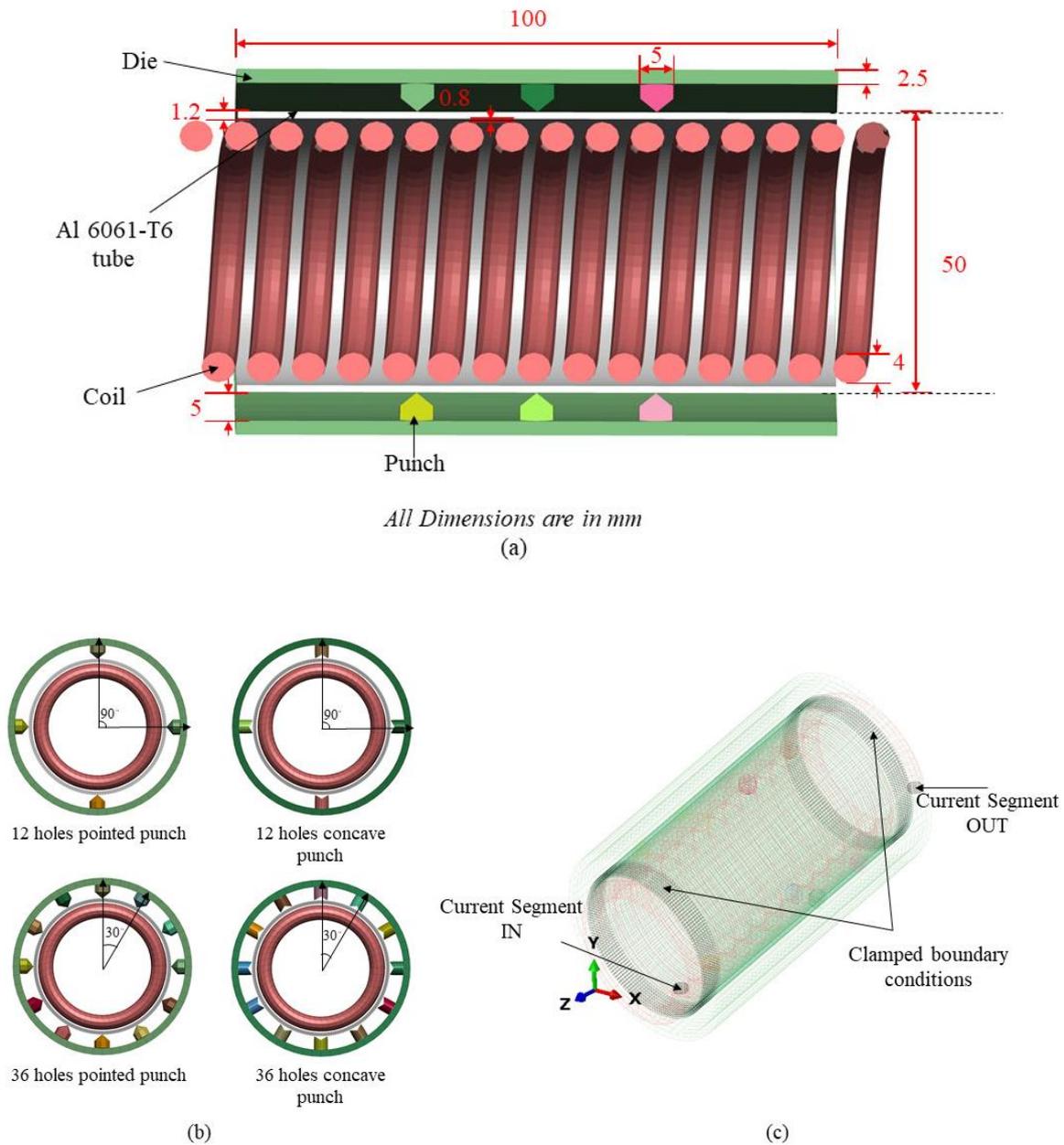

**Figure 4.** Pictorial representations of (a) geometrical dimensions of setup, (b) arrangement of the type of punches utilized, and (c) wireframe mode to locate input and boundary conditions.

In numerical simulations, the tube is made up of Al6061-T6 material, Copper solenoid coil, and SS304 punch and die arrangement. The material properties of the Al6061-T6 tube, Cu coil, and SS304 are illustrated in Table 2. To capture the deformations behavior for the high strain rate process, Johnson-Cook (JC) material model is used [26]. The empirical relationship of JC model is expressed as

$$\sigma = \left(A + B\varepsilon^n\right)\left[1 + C\ln\left(\frac{\dot{\varepsilon}}{\dot{\varepsilon}_0}\right)\right]\left[1 - \left(\frac{T - T_0}{T_m - T_0}\right)^m\right], \tag{14}$$

where $\sigma$ is equivalent stress (MPa), $\varepsilon$ is an equivalent plastic strain (mm/mm), $\dot{\varepsilon}$ and $\dot{\varepsilon}_0$ are the equivalent plastic strain rate and reference equivalent plastic strain rate (s$^{-1}$), respectively. *A, B, C, n*, and *m* are material constants, and their values are illustrated in Table 3. Also, *T, T$_0$*, and *T$_m$* represent the temperature, room temperature, and melting temperature (°C), respectively. The Johnson-cook failure model is used to capture the perforations of the tube. The fracture strain ($\varepsilon_f$) is sensitive to strain rate ($\dot{\varepsilon}$), stress triaxiality ratio ($\sigma^* = \sigma_m/\sigma_Y$), absolute temperature ($T^*$), and damage constants ($D_1$ to $D_5$) which are obtained from the experiment and can be expressed as [27]

$$\varepsilon^f = \left[D_1 + D_2 \exp\left(D_3 \sigma^*\right)\right]\left[1 + D_4 \ln\left(\dot{\varepsilon}^*\right)\right]\left[1 + D_5 T^*\right]. \tag{15}$$

*D* is a damage parameter that ranges from 0 (no damage) to 1 (complete damage and hence material deletion). The damage value is accumulated as expressed by Eq. (16),

$$D = \sum \frac{\Delta \varepsilon^p}{\varepsilon^f}, \tag{16}$$

where $\Delta \varepsilon^p$ is incremental equivalent plastic strain and $\varepsilon^f$ is the fracture strain. The damage constants values are depicted in Table 4. In the electromagnetic forming and perforation (EMFP) of the Al6061-T6 tube, the current waveform data was collected using an oscilloscope and a Rogowski coil for three different discharge energy values, as shown in Figure 5. The input current obtained from the experimental setup was assigned in the simulation model at the two extreme faces of the solenoid coil (as seen in Figure 4(c)). Since, the initial half pulse of the sinusoidal damped current curve was considered, as it potentially produces the effective deformation and perforation [28], it is sufficient to simulate up to 30 μs of current as shown in Figure 5. A timestep of 0.1 μs was used during simulations, with a maximum number of iterations of 1000, and the total CPU running time was 7 hours.

**Table 2.** Properties of material used.

| Properties | SI Units | Al 6061-T6 | Cu | SS304 |
|---|---|---|---|---|
| Density ($\rho$) | kg/m$^3$ | 2703 | 8940 | 7900 |
| Young's Modulus (E) | GPa | 69.8 | 97 | 210 |
| Shear Modulus (G) | GPa | 26 | 45 | 70 |
| Poisson's ratio ($\upsilon$) | | 0.33 | 0.33 | 0.3 |
| Electrical conductivity ($1/\gamma$) | MS/m | 25 | 58 | 1.4 |
| Thermal conductivity ($k_T$) | W/m-K | 126 | 391 | 14 |

**Table 3.** Johnson-Cook (JC) plasticity model: input constants. [26]

| Material | A [MPa] | B [MPa] | C | n | m | $T_m$ (K) | $\dot{\varepsilon}_0$ [s$^{-1}$] |
|---|---|---|---|---|---|---|---|
| Al 6061-T6 | 324.1 | 113.8 | 0.002 | 0.42 | 1.34 | 925 | 1.0 |

**Table 4.** Johnson-Cook (JC) dynamic failure model: damage constants. [27]

| Material | $D_1$ | $D_2$ | $D_3$ | $D_4$ | $D_5$ |
|---|---|---|---|---|---|
| Al 6061-T6 | −0.77 | 1.45 | −0.47 | 0.0 | 1.6 |

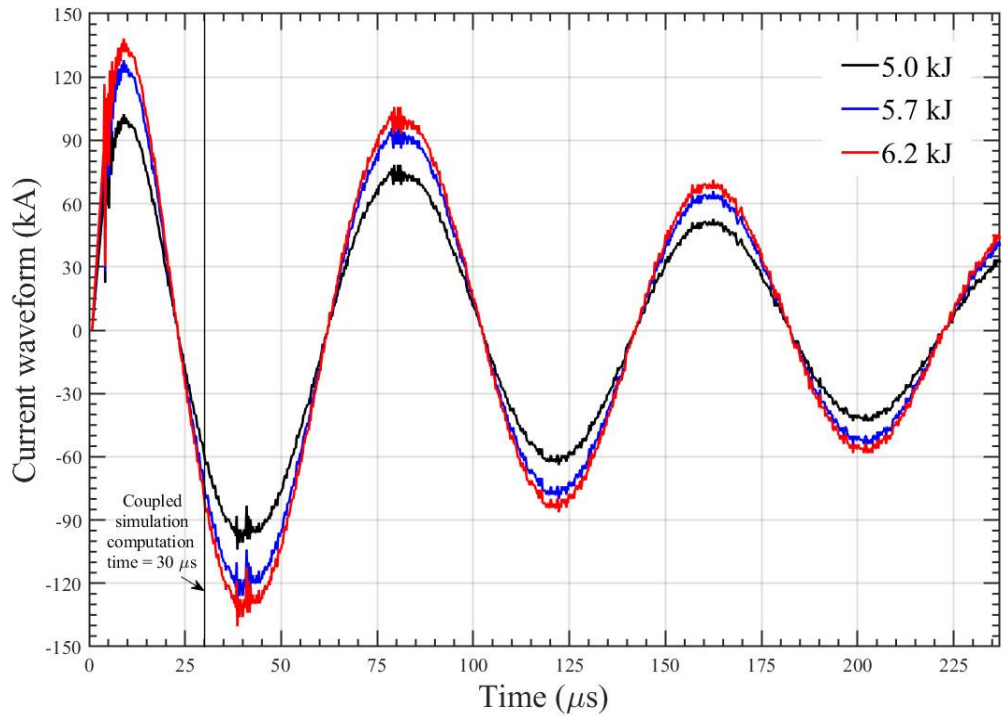

**Figure 5.** Instrumented measured current vs time.

## 4. Results and discussions:

Experiments were conducted out on the Al6061-T6 tube using pointed and concave punches for energy levels varied from 4.8 to 6.5 kJ [5]. At lower energy levels, there is no perforation and partial indentations at the punch locations. The tube is entirely perforated by increasing the discharge energy, and the required hole diameter is achieved. During perforation, pointed punches exhibit petaling phenomena whereas blank or slug separation is observed with concave punches [25].

*4.1 Electromagnetic and Deformation analysis:*

As electromagnetic forming and perforations occur in a few microseconds, capturing the variations in process variables, including Lorentz forces, resultant velocity, equivalent stress, and effective plastic strain during experiments, is not possible. In order to shed light on EMFP physics, the coupled simulation analysis is used for discharge energy values of 5.7 kJ. We have selected four different elements surrounding the perforated hole. Various process variables have been presented as an average over these elements. Effect of multiple punches, different punch types and axial location of the holes have been thoroughly analysed in the following sections.

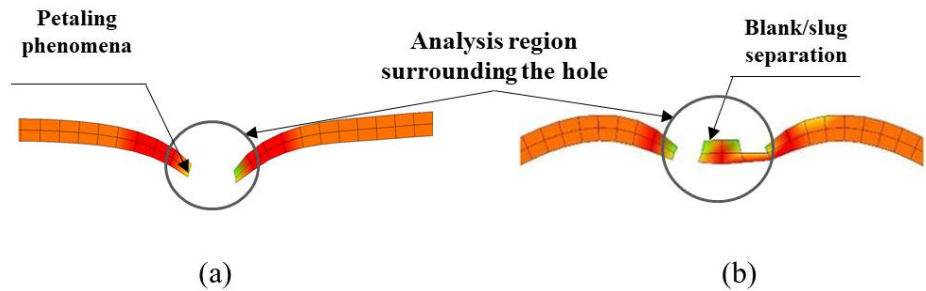

**Figure 6.** Selected elements in (a) pointed punch and (b) concave punch simulations, to analyse the variation of different parameters.

*4.1.1. Lorentz force distributions based on type of punches and number of punches:*

Type of punches can influence the temporal variation of the Lorentz force as well as the timepoint of the perforation initiation. For example, a pointed punch generates a higher concentration of Lorentz force at its tip, while a concave punch spreads the Lorentz force over a larger area of the tube. Increasing the number of the perforated holes can result in a reduction in Lorentz force intensity for a specified energy level. This is because the current flow is divided among multiple holes, reducing the concentration of current and the resultant Lorentz force. For an instance, the average Lorentz force distributions with respect to time for a discharge energy value of 5.7 kJ for four specific cases are depicted in Figure 7, which validate the above theories. It is also observed that the distribution of this force is generally higher at the center of the tube than at the ends in accordance with the electromagnetic field distribution. As the number of punches increases, effective force is shared by more indentation locations.

Hence, the average value of Lorentz force reduces.

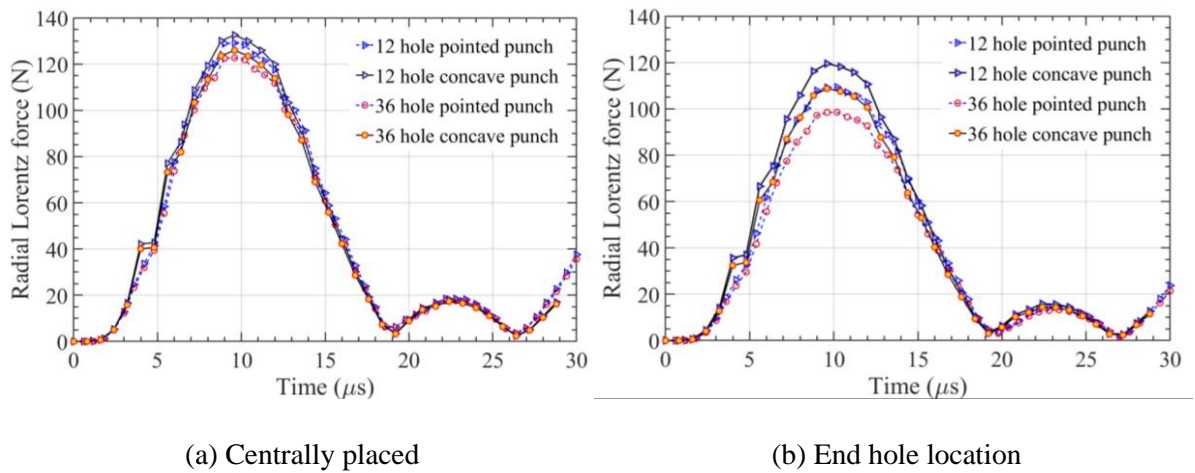

(a) Centrally placed           (b) End hole location

**Figure 7.** Variations of average Lorentz force distributions with respect to time for 5.7 kJ discharge energy for four different cases.

*4.1.2. Resultant velocity distribution based on type of punches and number of punches:*

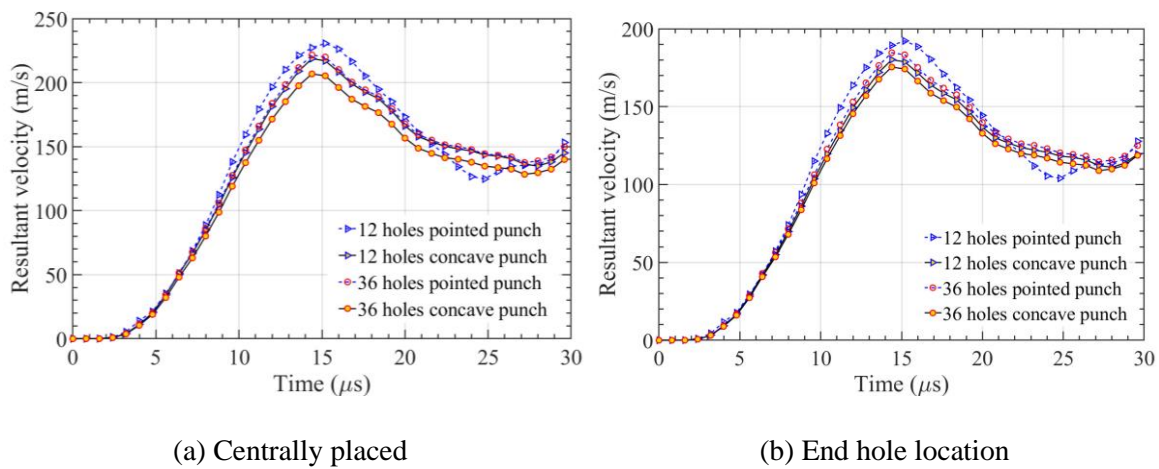

(a) Centrally placed           (b) End hole location

**Figure 8.** Variations of average resultant velocity with respect to time for 5.7 kJ discharge energy for four different cases.

The results show that increasing the number of punches in the electromagnetic forming and perforation process leads to a reduction in the magnitude of the resultant velocity. This can be attributed to the spreading of energy over a larger area as the number of punches increases. The simulation findings also show that the use of concave punches results in a slightly lower resultant velocity as compared to the use of pointed punches due to the rise in the shearing area. The distribution of Lorentz force is generally higher at the center of the tube than at the ends as discussed in previous section. This higher force results in a higher velocity at the center of the tube as compared to the ends.

*4.1.3. von-Mises stress and effective plastic strain variations based on type of punches and number of punches:*

A pointed punch generates higher stresses than a concave punch due to more concentrated load distribution. Also, a concave tool generates more uniform stresses, leading to improved material flow and formability. The variations in equivalent stress (von Mises stress) and effective plastic strain along time are depicted in Figure 9 and Figure 10 respectively. With increase in number of punches from 12 to 36, effective Lorentz force in one particular punch location reduces. This results in reduction in maximum von-Mises stress and effective plastic strain. This is observed for both pointed and concave punches. Also, we have more effective plastic strain and associated von-Mises stress in the center location than those in the end location. This is due to the higher Lorentz forces in the central region.

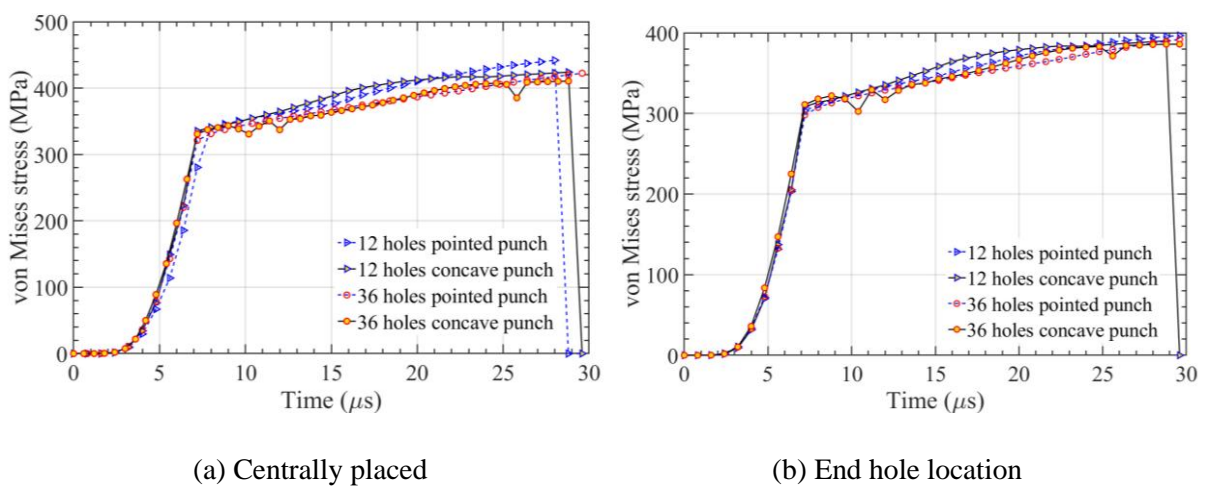

(a) Centrally placed  (b) End hole location

**Figure 9.** Variations of average equivalent stress with respect to time for 5.7 kJ energy level for four distinct cases.

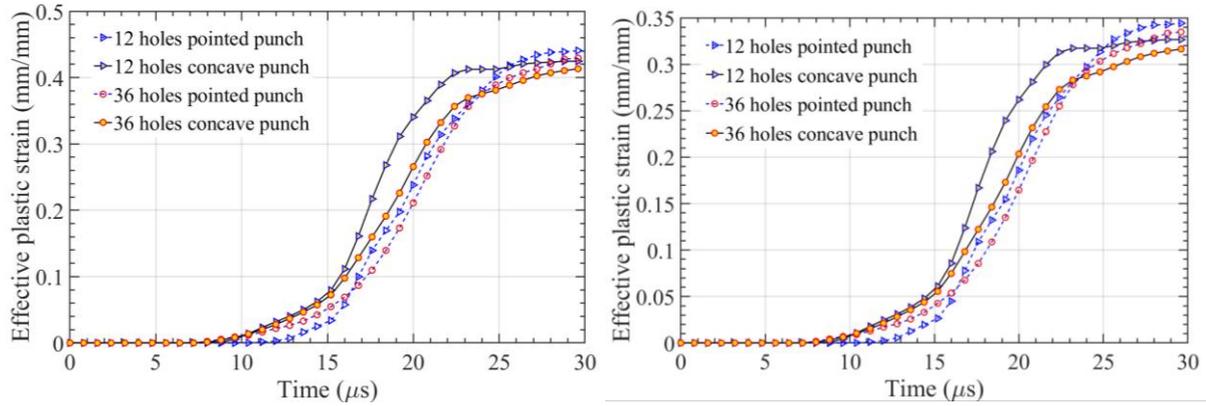

(a) Centrally placed  (b) End hole location

**Figure 10.** Variations of average effective plastic strain with respect to time for 5.7 kJ energy level for four distinct cases.

The effective plastic strain and von Mises stress values, obtained through numerical simulations, provide important insight about the material behaviour in the process of forming and perforation.

*4.2 Validation and Comparison of coupled FE model results with experimental results:*

*4.2.1. Comparison of average hole diameter:*

The Ls-Dyna™ results and the experimental results [5] are compared for multi-point EMFP for average hole diameter for three discharge energies as illustrated in Figure 11. As the energy level increases, the amount of magnetic force on the tube also increases, leading to an increase in the average diameter of the hole. The average diameter of the hole produced by a concave punch is typically larger than that produced by a pointed punch for same discharge energy. A concave punch has a curved surface that spreads the force over a larger shearing area, leading to a larger hole diameter. A pointed punch, on the other hand, has a sharp, pointed tip that concentrates the force of the punch in a smaller shearing area, leading to a smaller hole diameter. When a tube is perforated with multiple holes, the effective magnetic Lorentz force is shared by multiple holes. This results in a reduction in the plastic deformation behavior of the material around each hole, which results in reduction in the final hole diameter. Additionally, the presence of neighbouring holes also affects the punching process and the final hole diameter. The close alignment among numerical and experimental results shows the efficacy of the numerical simulations. The diameter of the holes increases with the rise in energy values. However, there is a threshold limit beyond which the workpiece can fracture between two adjacent holes in the axial direction. The threshold limit is 6.6 kJ for 12 holes perforated tube and 6.2 kJ for 36 holes perforated tube with a workpiece thickness of 1.2 mm [5].

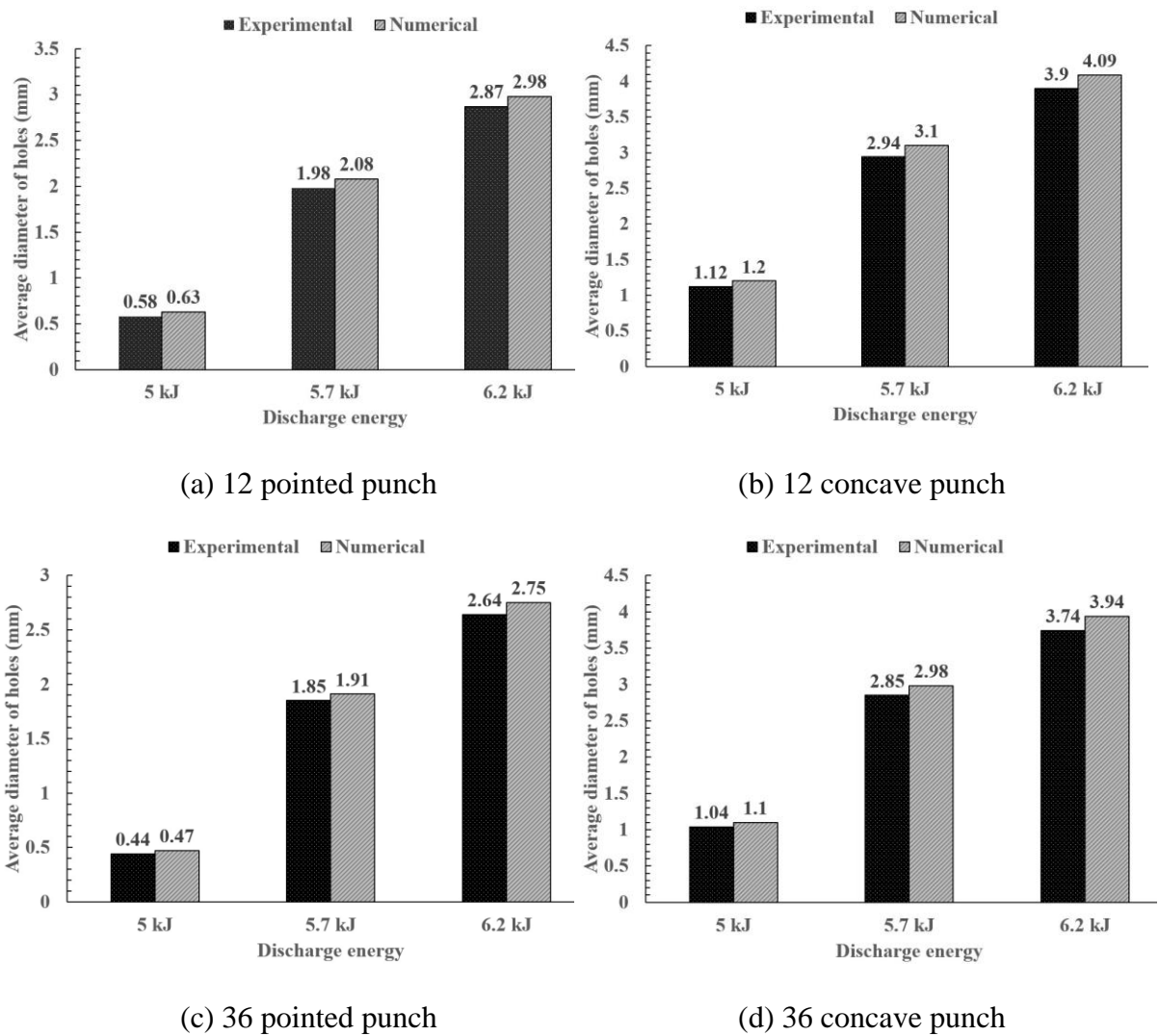

(a) 12 pointed punch    (b) 12 concave punch

(c) 36 pointed punch    (d) 36 concave punch

**Figure 11.** Comparison of experimental and numerical results for the average diameter of holes

*4.2.2. Comparative study on number of perforated holes:*

Number of perforated holes for different discharge energies have been compared between numerical simulation and experiments [5] in Fig. 12. Higher discharge energy leads to the formation of more holes, as the intensity of the perforating force increases throughout the tube surfaces. A pointed punch with a small shearing area, perforates more easily than a concave punch with a larger shearing area as the pointed punch concentrates the force in a smaller area and reduces the material flow into the shearing area. A concave punch requires more energy for blank separation as compared to that for a pointed punch, because the concave punch removes larger amount of material due to its larger shearing area. In the case of multiple holes, the difficulty of blank separation is compounded as the presence of neighbouring holes affect the ability to cleanly separate the blanks. Therefore, for more no. of punches (Here 36) at the same discharge energies, numbers of complete holes in case of pointed punches are considerably higher than that in case of concave punches. Again, our numerical model is highly capable to simulate the experimental results of number of perforated holes for different discharge energies and for both the punches.

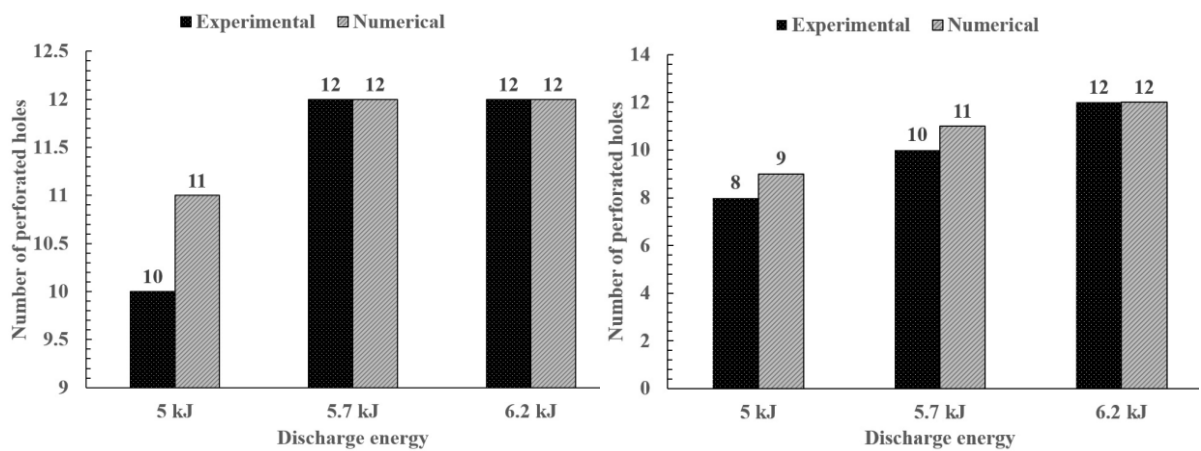

(a) 12 pointed punch          (b) 12 concave punch

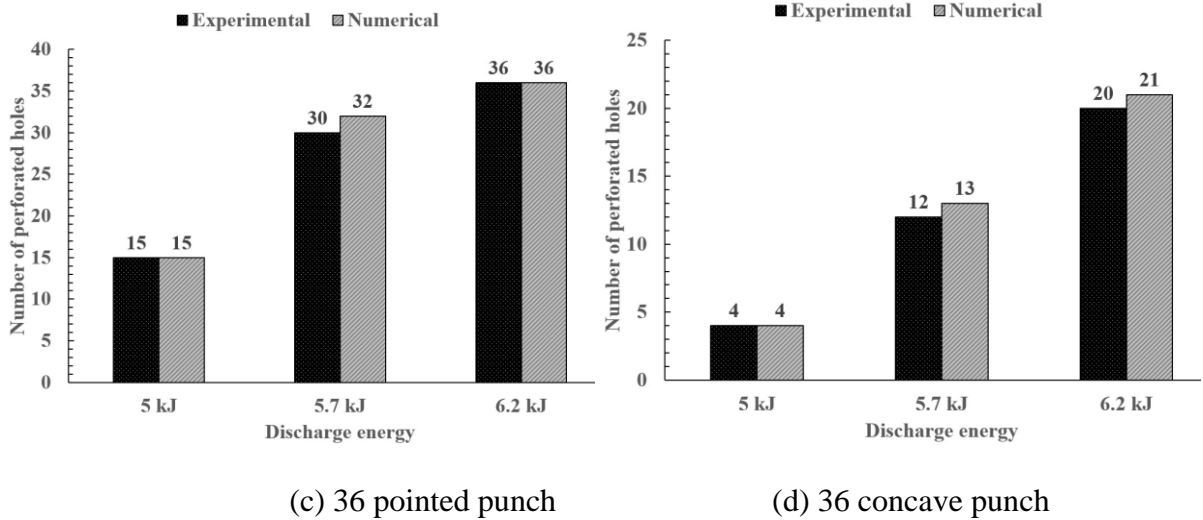

(c) 36 pointed punch    (d) 36 concave punch

**Figure 12.** Comparison of experimental and numerical results for a number of perforated holes.

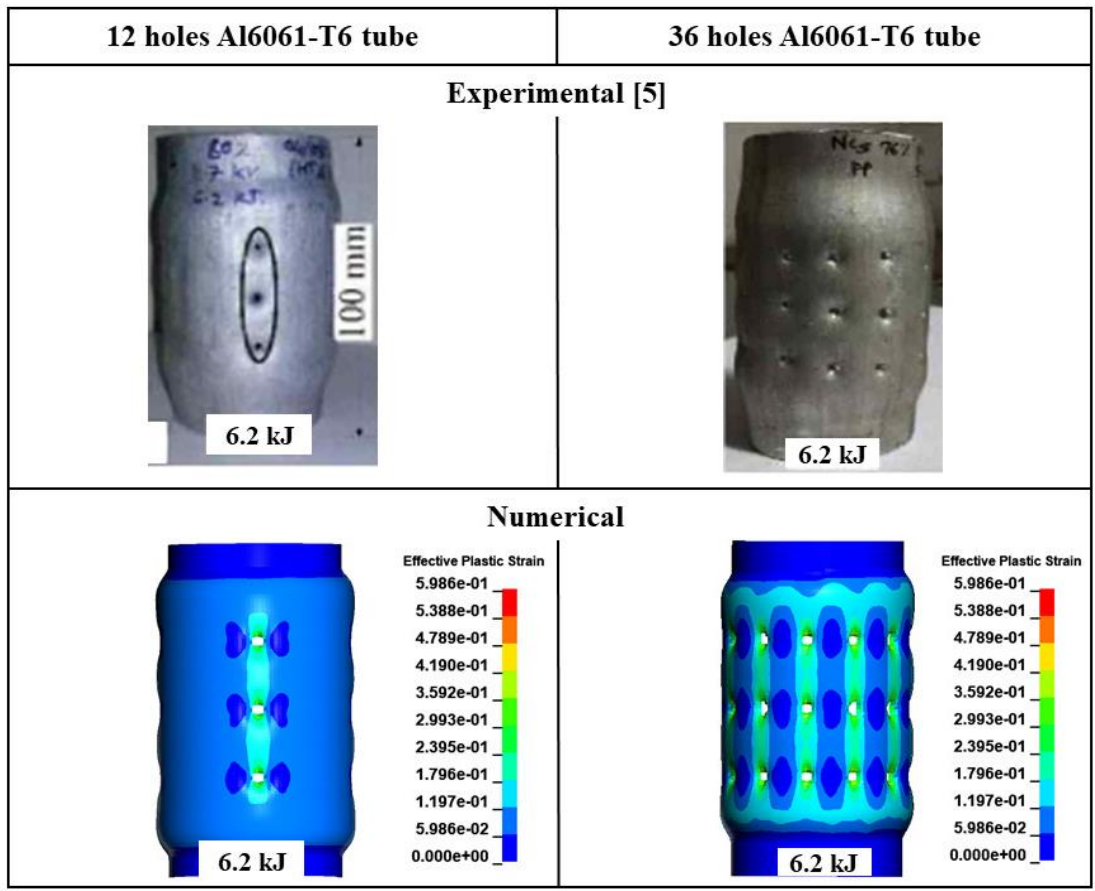

**Figure 13.** Final profile of deformed and perforated tube in experiments and simulations for pointed punch.

Figure 13 represent the deformation and perforation profile of the tubes under the action of multiple pointed punches, where again close agreement between experiments and simulations is evident. The results indicate that despite using higher discharge energy values, there is sign of petal phenomena and smaller hole diameters without any blank removal.

Figure 14 represent the close agreement between experimental and numerical deformed profiles after EMFP of tubes with 12 and 36 concave punches at 5.7 kJ and 6.2 kJ discharge energies respectively. In summary, the results indicate that multiple concave punches lead to separation of the blank or slug in the sheared zone but require higher discharge energy to achieve full removal. There is also an increase in hole diameter as compared to the use of pointed punches at same discharge energies.

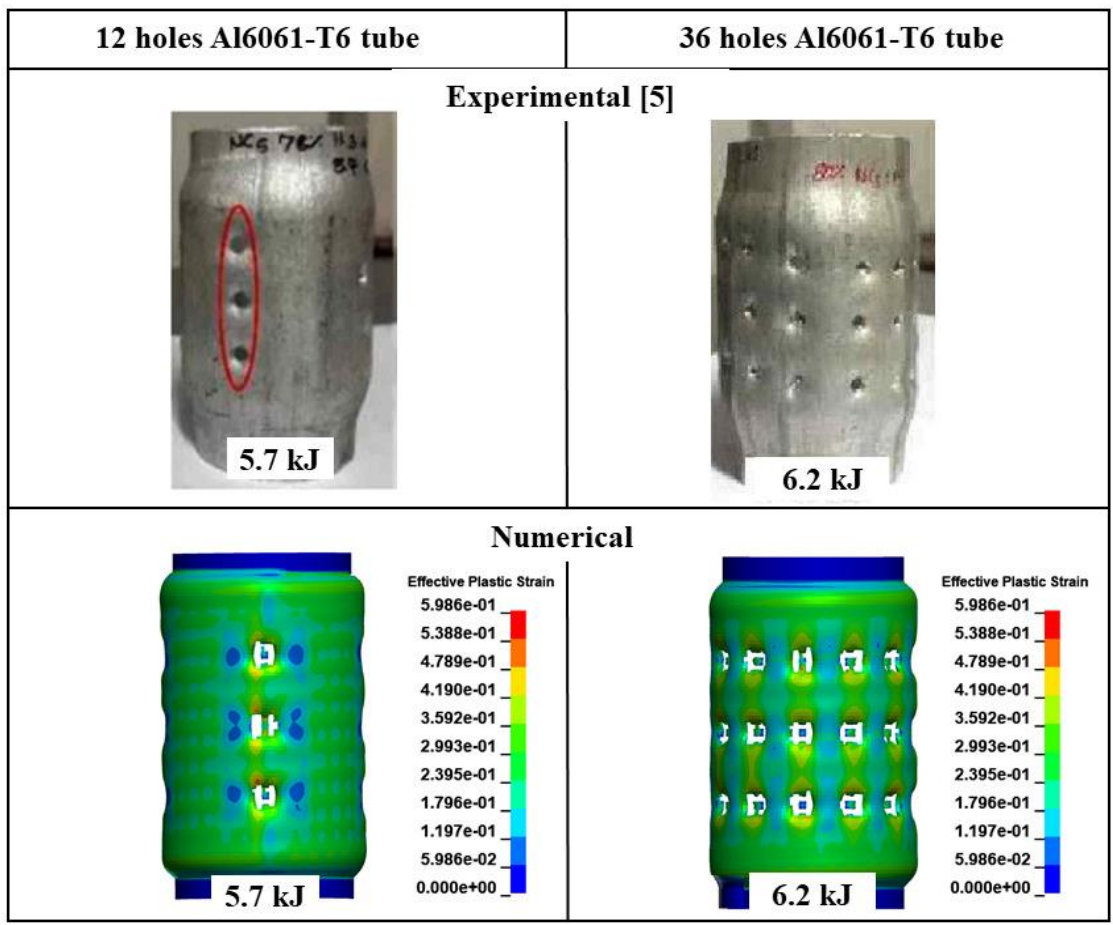

**Figure 14.** Final profile of deformed and perforated tube in the experiments and simulations for concave punch.

**Conclusions:**

A coupled finite element is used in this study to analyze the multi-physics of the EMFP of Al6061-T6 tubes. For a 12 and 36 holes pointed and concave punch design, the critical output process parameters, such as Lorentz force distributions and deformation behavior during the process, are explored. The numerical simulations are conducted out for three distinct discharge energy values to understand the variations of hole diameters and no. of complete perforated holes which are essential quality measures for EMFP of Al 6061-T6 tubes. Following are the brief conclusions of the above analysis.

- Among all the numerical simulations with either pointed or concave punches, and for creating both 12 and 36 perforated holes, there is very firm agreement among Ls-Dyna™ and experimental results with maximum deviation of less than 6%.
- For a given discharge energy, the shape of the punch can have a significant consequence on the amplitude and distribution of the Lorentz force and, hence on, the initiation of perforation. A pointed punch generates a higher of Lorentz force as compared to that by a concave punch. It is due to the larger shearing area with the concave punches. Lorentz force also reduces with the increase of the number of punches due to the increase in the total shearing area.
- Due to the higher Lorentz forces, we have higher effective plastic strain and higher von-Mises stress for less no. of punches, and for pointed punch as opposed to the concave counterpart.
- A pointed punch results in a higher resultant velocity as compared to that of a concave punch due to its smaller surface area. The higher velocity from a pointed punch results in rapid initiation of perforation but resulted in smaller diameter of holes. On the other hand, a concave punch provides a larger surface area, resulting in a lower velocity and a slower initiation of perforation, and shallower penetrations causing larger hole diameter.

- The manufacturers can obtain a guideline from numerical simulations to maintain a balance between the average diameter of holes and the number of perforated for a given tube and punch geometries and discharge energy. It will save a large amount of resources eliminating the repeated experiments.

**Disclosure statement**

No potential conflict of interest between the authors.